\documentclass[letterpaper]{article}
\usepackage[T1]{fontenc}
\usepackage{lmodern}
\usepackage{amsmath,amsfonts,amsthm}
\usepackage[margin=1in]{geometry}
\usepackage[
colorlinks=true,
linkcolor=blue,
citecolor=blue,
urlcolor=blue
]{hyperref}
\newtheorem{theorem}{Theorem}[section]

\newtheorem{corollary}[theorem]{Corollary}

\newcommand{\C}{\mathbb C}
\newcommand{\Rcal}{\mathcal R}
\newcommand{\norm}[2]{\left\lVert #1\right\rVert_{#2}}
\DeclareMathOperator{\dist}{dist}
\DeclareMathOperator{\diam}{diam}
\newcommand{\capE}{\operatorname{cap}(E)}

\title{Product Asymptotics for Rows of the Walsh Table on Regular Compact Sets}
\author{V. A. Prokhorov}
\date{}

\begin{document}
	
	\maketitle
	
	\begingroup
	\renewcommand{\thefootnote}{}
	
	\footnotetext{\textit{2020 Mathematics Subject Classification.}
		Primary 41A20, 41A21; Secondary 30E10.}
	
	\footnotetext{\textit{Key words and phrases.}
		Best rational approximation. Rows of the Walsh table. Hadamard-type theorems
		for meromorphy radii. Regular compact sets.}
	\endgroup
	
	\begin{abstract}
		Let $E$ be a regular compact set with connected
		complement, and let $f$ be holomorphic
		in a neighborhood of $E$.
		We establish Hadamard-type asymptotic formulas
		for products of errors of best uniform rational
		approximation of $f$ on $E$, formed from
		the first $m+1$ rows of the Walsh table. When the corresponding meromorphy radii
		are finite, we also show that the extremal
		exponential rates of the approximation errors
		along these rows are attained
		along a common subsequence. 
		We also prove a sharp asymptotic comparison
		theorem for ratios of best rational approximation
		errors along fixed rows of the Walsh table
		on $E$ and its closed Green sublevel sets,
		and derive corresponding sharp asymptotic
		estimates for ratios of their products.
		These results extend earlier product theorems from continua with connected complement and Jordan boundary to the more general setting of regular compact sets with connected complement.
	\end{abstract}

\section{Notation}

Let $E\subset\mathbb{C}$ be a compact set
of positive logarithmic capacity such that
$\Omega:=\widehat{\mathbb{C}}\setminus E$ is connected
and regular with respect to the Dirichlet problem.
We denote by $g_\Omega(z,\infty)$ the Green function
of $\Omega$ with pole at infinity, extended
continuously to $\mathbb{C}$ by zero on $E$.

For each $R>1$, define the open Green sublevel set
\[
G_R:=\{z\in\C:
g_\Omega(z,\infty)<\log R\}, 
\]
and the corresponding closed Green sublevel set
\[
E_R:=\overline{G_R}
=\{z\in\C:
g_\Omega(z,\infty)\le\log R\}.
\] 
We also set $E_1:=E$.

For nonnegative integers $n$ and $m$, let
$\Rcal_{n,m}$ denote the class of rational
functions of type $(n,m)$, that is,
\[
\Rcal_{n,m}:=
\left\{\frac pq:
\deg p\le n,\ \deg q\le m,\
q\not\equiv0\right\}.
\]
For a nonempty compact set $K\subset\mathbb{C}$
and a function $g$ continuous on $K$, define
\[
\rho_{n,m}(g;K):=
\inf_{r\in\Rcal_{n,m}}\norm{g-r}{K},
\]
the error of best uniform approximation of $g$ on $K$
by rational functions from the class $\Rcal_{n,m}$, where
\[
\norm{h}{K}:=\max_{z\in K}|h(z)|.
\]

Now let $f$ be holomorphic in a neighborhood
of $E$. For each integer $m\ge0$, define
$R_m=R_m(f;E)$ as the supremum of all $R>1$
for which $f$ admits a meromorphic continuation
to $G_R$ having at most $m$ poles in total,
counted with multiplicity across all
components of $G_R$. In particular, when $m=0$,
the continuation is holomorphic.
The radii of meromorphic continuation $R_m$
are nondecreasing in $m$, that is,
\[
1<R_0\le R_1\le R_2\le\cdots.
\]

\section{Classical results and motivation}

For each integer $k\geq 0$, the sequence
$\{\rho_{n,k}(f;E)\}_{n\geq 0}$ is called
the $k$th row of the Walsh table.

The classical results of Saff~\cite{S1971}
and Gonchar~\cite{G1974} establish a fundamental
connection between the asymptotic behavior
of the rows of the Walsh table and meromorphic
continuation. Saff established this connection
for continua with connected complement and Jordan boundary,
and Gonchar subsequently extended it to regular
compact sets with connected complement.
The Saff--Gonchar theorem asserts that, for
every fixed $k\geq 0$,
\[
\limsup_{n\to\infty}
\rho_{n,k}(f;E)^{1/n}
=\frac{1}{R_k}.
\]
Thus, the radius $R_k$ of meromorphic continuation
determines, and is determined by, the asymptotic
exponential rate of best rational approximation
along the $k$th row of the Walsh table.
For $k=0$, this reduces to the classical
Bernstein--Walsh theorem, which characterizes
analytic continuation in terms of the rate
of best polynomial approximation.

Another classical result relevant to our study
is Hadamard's theorem~\cite{H1892}, which expresses
the radii of meromorphic continuation of a
power series through the asymptotic behavior
of its Hankel determinants.

Let $F$ be holomorphic in a neighborhood of
the origin, with Taylor expansion
\[
F(z)=\sum_{j=0}^{\infty}F_jz^j.
\]
Denote by $r_0>0$ the radius of convergence
of this series. For each integer $m\geq 1$,
let $r_m$ be the supremum of all $r>0$ for which
$F$ admits a meromorphic continuation to
the disk $|z|<r$ with at most $m$ poles,
counted with multiplicity.

For integers $m\geq 1$ and $n\geq 0$, consider
the Hankel matrix
\[
K_{n,m}:=
\begin{bmatrix}
	F_{n+m-1} & F_{n+m-2} & \cdots & F_n\\
	F_{n+m-2} & F_{n+m-3} & \cdots & F_{n-1}\\
	\vdots & \vdots & \ddots & \vdots\\
	F_n & F_{n-1} & \cdots & F_{n-m+1}
\end{bmatrix},
\]
where we set $F_j:=0$ for $j<0$.
Hadamard's theorem states that
\[
\limsup_{n\to\infty}
|\det K_{n,m+1}|^{1/n}
=\frac{1}{r_0r_1\cdots r_m},
\qquad m\geq 0.
\]
The exponential asymptotics of these determinants
depend on the product of the first $m+1$
meromorphy radii.

This leads to the question of whether the product
$R_0R_1\cdots R_m$ also determines the exponential
asymptotics of products of best rational approximation
errors associated with the first $m+1$ rows
of the Walsh table.

To investigate this issue, for $n\geq m$
we consider the product
\[
\prod_{k=0}^{m}
\rho_{n-m+k,k}(f;E),
\]
formed by taking one approximation error
from each of the first $m+1$ rows of the Walsh table,
with the first index increasing from $n-m$ to $n$.
Applying the Saff--Gonchar theorem to the
individual factors gives
\begin{equation}\label{L1}
	\limsup_{n\to\infty}
	\left(
	\prod_{k=0}^{m}
	\rho_{n-m+k,k}(f;E)
	\right)^{1/n}
	\leq \frac{1}{R_0R_1\cdots R_m}.
\end{equation}
The analogy with Hadamard's determinant formula
suggests that equality should hold in~\eqref{L1}.
The reverse inequality, however, does not
follow from the individual-row asymptotics alone,
since they do not guarantee that the extremal
rates are attained along a common subsequence.

In~\cite{Pro2026}, Hadamard-type asymptotics
for products of best rational approximation
errors were obtained for continua with
connected complement and Jordan boundary.
The analysis combined the classical theorems
of Hadamard and Saff--Gonchar
with weighted Hankel operators and a weighted
Adamyan--Arov--Kre\u{\i}n-type theorem
for meromorphic approximation.
Products of singular numbers of weighted Hankel
operators were represented by determinants
associated with bilinear Hankel forms,
and contour deformation played a central role
in obtaining comparison estimates.
The results also included sharp asymptotic
estimates for ratios of the corresponding
products and the existence of a common
extremal subsequence.

The main result of the present paper establishes
the Hadamard-type asymptotic formula for products
of best rational approximation errors on regular
compact sets with connected complement.
In particular, $E$ need not be connected,
and no additional geometric conditions
are imposed on its boundary. Unlike the operator-theoretic
approach of~\cite{Pro2026}, our proofs
work directly with best rational approximation
errors.

For the Hadamard-type product formula,
the essential ingredient is the existence
of ordinary limits for the $n$th roots
of the approximation errors along the rows
with $R_k<R_m$.
Our proof of these limits is based on
ideas developed by Saff~\cite{Saff1969a,Saff1969b}.
In particular, we adapt his method of
comparing rational approximants, clearing
denominators, applying the Bernstein--Walsh
inequality, and using a pole-counting argument
to obtain a contradiction.
We compare best rational approximants
from different rows of the Walsh table
and use Gonchar's analogue of the
de Montessus de Ballore theorem
to control the convergence of their
denominators.
Together with Rouch\'e's theorem,
these arguments yield the ordinary limits
required for a direct proof of the product formula. 
When $R_m<\infty$, the existence of a common
extremal subsequence for the first $m+1$
rows of the Walsh table follows directly
from the exact product asymptotics.

Questions concerning ratios of best rational
approximation errors were previously studied
by Prokhorov~\cite{Pro2005}
and by Kouchekian and
Prokhorov~\cite{KouchPro2009}.
In the present paper, we also establish
sharp asymptotic estimates for ratios
of individual approximation errors
on $E$ and $E_R$ for regular compact sets
with connected complement.
These estimates yield corresponding
bounds for ratios of their products.
The comparison of approximation errors on $E$
and its closed Green sublevel sets $E_R$
is obtained by polynomial interpolation
at Fekete points.

\section{Product asymptotics for regular compact sets}

In this section, we present our results.
Our main theorem shows that the upper estimate
in~\eqref{L1} is sharp.

\begin{theorem}\label{L2}
	Let $E$ be a regular compact set with connected
	complement, and let $f$ be holomorphic
	in a neighborhood of $E$. Fix $m\geq 0$. Then
	\begin{equation}\label{L3}
		\limsup_{n\to\infty}
		\left(
		\prod_{k=0}^{m}
		\rho_{n-m+k,k}(f;E)
		\right)^{1/n}
		=
		\frac{1}{R_0R_1\cdots R_m}.
	\end{equation}
\end{theorem}
For $1<R<R_0$, the Green function associated
with the complement of $E_R$ satisfies
\[
g_{\widehat{\mathbb C}\setminus E_R}(z,\infty)
=
g_\Omega(z,\infty)-\log R,
\qquad z\in\mathbb C\setminus E_R.
\]
Consequently, the radii of meromorphic
continuation of $f$ relative to $E_R$
are $R_k/R$, $k\geq 0$.
By the Saff--Gonchar theorem, for each fixed
$k\geq 0$,
\begin{equation}\label{L4}
	\limsup_{n\to\infty}
	\rho_{n,k}(f;E_R)^{1/n}
	=
	\frac{R}{R_k}.
\end{equation}
Applying Theorem~\ref{L2}
to $E_R$ yields the following corollary.

\begin{corollary}
	Let $E$ be a regular compact set with connected
	complement, and let $f$ be holomorphic
	in a neighborhood of $E$.
	Fix $m\geq 0$. Then, for every
	$1\leq R<R_0$,
	\begin{equation}\label{L5}
		\limsup_{n\to\infty}
		\left(
		\prod_{k=0}^{m}
		\rho_{n-m+k,k}(f;E_R)
		\right)^{1/n}
		=
		\frac{R^{m+1}}{R_0R_1\cdots R_m}.
	\end{equation}
\end{corollary}

We next consider a comparison theorem for best
rational approximation errors along a fixed
row of the Walsh table. Since $E\subset E_R$,
we have
\[
\rho_{n,m}(f;E)\leq \rho_{n,m}(f;E_R).
\]

The following result provides a sharper asymptotic
estimate for the ratio of these errors and
shows that equality holds when $R_m<\infty$.

\begin{theorem}\label{L6}
	Let $E$ be a regular compact set with connected
	complement, and let $f$ be holomorphic
	in a neighborhood of $E$.
	Fix $m\geq 0$ and $1<R<R_0$. 
	Assume that $f$ is not
	a rational function with at most $m$
	finite poles, counted with multiplicity. Then
	\begin{equation}\label{L7}
		\limsup_{n\to\infty}
		\left(
		\frac{\rho_{n,m}(f;E)}
		{\rho_{n,m}(f;E_R)}
		\right)^{1/n}
		\leq \frac{1}{R}.
	\end{equation}
	Moreover, if $R_m<\infty$, then
	\begin{equation}\label{L8}
		\limsup_{n\to\infty}
		\left(
		\frac{\rho_{n,m}(f;E)}
		{\rho_{n,m}(f;E_R)}
		\right)^{1/n}
		=\frac{1}{R}.
	\end{equation}
\end{theorem}

To see the equality in~\eqref{L8}, suppose that
$R_m<\infty$. We use the elementary fact
that, for two sequences of positive numbers
$\{a_n\}$ and $\{b_n\}$ with finite upper limits
and $\limsup_{n\to\infty}b_n>0$, one has
\[
\limsup_{n\to\infty}\frac{a_n}{b_n}
\geq
\frac{\limsup_{n\to\infty}a_n}
{\limsup_{n\to\infty}b_n}.
\]
Applying the Saff--Gonchar theorem on $E$
and $E_R$ (see~\eqref{L4}), we obtain
\[
\limsup_{n\to\infty}
\left(
\frac{\rho_{n,m}(f;E)}
{\rho_{n,m}(f;E_R)}
\right)^{1/n}
\geq
\frac{1/R_m}{R/R_m}
=\frac{1}{R}.
\]
Together with \eqref{L7}, this proves
\eqref{L8}.

The fixed-row comparison also yields an estimate
for products involving several rows. Indeed,
applying Theorem~\ref{L6}
with $m$ replaced by $k$, for $k=0,\ldots,m$,
and with numerator degree $n-m+k$, gives
\[
\limsup_{n\to\infty}
\left(
\frac{\rho_{n-m+k,k}(f;E)}
{\rho_{n-m+k,k}(f;E_R)}
\right)^{1/n}
\leq \frac{1}{R}.
\]
Multiplying these inequalities leads to
the following product estimate.

\begin{corollary}
	Let $E$ be a regular compact set with connected
	complement, and let $f$ be holomorphic
	in a neighborhood of $E$.
	Fix $m\geq 0$ and $1<R<R_0$. 
	Assume that $f$ is not
	a rational function with at most $m$
	finite poles, counted with multiplicity. Then
	\begin{equation}\label{L9}
		\limsup_{n\to\infty}
		\left(
		\prod_{k=0}^{m}
		\frac{\rho_{n-m+k,k}(f;E)}
		{\rho_{n-m+k,k}(f;E_R)}
		\right)^{1/n}
		\leq \frac{1}{R^{m+1}}.
	\end{equation}
	Moreover, if $R_m<\infty$, then
	\begin{equation}\label{L10}
		\limsup_{n\to\infty}
		\left(
		\prod_{k=0}^{m}
		\frac{\rho_{n-m+k,k}(f;E)}
		{\rho_{n-m+k,k}(f;E_R)}
		\right)^{1/n}
		=
		\frac{1}{R^{m+1}}.
	\end{equation}
\end{corollary}
The equality in~\eqref{L10} follows
from the same elementary inequality used above.
Indeed, when $R_m<\infty$, applying this
inequality to the $n$th roots of the
corresponding products on $E$ and $E_R$
and using~\eqref{L3} and~\eqref{L5},
we obtain

\[
\begin{aligned}
	&\limsup_{n\to\infty}
	\left(
	\prod_{k=0}^{m}
	\frac{\rho_{n-m+k,k}(f;E)}
	{\rho_{n-m+k,k}(f;E_R)}
	\right)^{1/n}\\
	&\qquad\geq
	\frac{1/(R_0R_1\cdots R_m)}
	{R^{m+1}/(R_0R_1\cdots R_m)}
	=\frac{1}{R^{m+1}}.
\end{aligned}
\]
Together with \eqref{L9}, this proves
\eqref{L10}.

In particular, when $R_m<\infty$, the upper limit
in~\eqref{L10} depends only on $R$ and $m$,
and is independent of the individual
meromorphy radii $R_0,\ldots,R_m$.

The preceding theorems describe asymptotic
behavior in terms of upper limits.
A further question is whether the extremal
exponential rates associated with the
different rows can be realized simultaneously.
The following theorem establishes the existence
of a common extremal subsequence for the
individual errors, their partial products,
and the corresponding quotients.

\begin{theorem}
	Let $E$ be a regular compact set with connected
	complement, and let $f$ be holomorphic
	in a neighborhood of $E$.
	Fix $m\geq 0$ and assume that $R_m<\infty$.
	Let $\Lambda\subset\mathbb N$ be a subsequence
	such that
	\[
	\lim_{\substack{n\to\infty\\n\in\Lambda}}
	\left(
	\prod_{j=0}^{m}
	\rho_{n-m+j,j}(f;E)
	\right)^{1/n}
	=
	\frac{1}{R_0R_1\cdots R_m}.
	\]
	Then, for every $k=0,1,\ldots,m$ and every
	$1\leq R<R_0$,
	\[
	\lim_{\substack{n\to\infty\\n\in\Lambda}}
	\rho_{n-m+k,k}(f;E_R)^{1/n}
	=
	\frac{R}{R_k},
	\]
	and
	\[
	\lim_{\substack{n\to\infty\\n\in\Lambda}}
	\left(
	\prod_{j=0}^{k}
	\rho_{n-m+j,j}(f;E_R)
	\right)^{1/n}
	=
	\frac{R^{k+1}}{R_0R_1\cdots R_k}.
	\]
	Moreover, for every $k=0,1,\ldots,m$
	and every $1<R<R_0$,
	\[
	\lim_{\substack{n\to\infty\\n\in\Lambda}}
	\left(
	\prod_{j=0}^{k}
	\frac{\rho_{n-m+j,j}(f;E)}
	{\rho_{n-m+j,j}(f;E_R)}
	\right)^{1/n}
	=
	\frac{1}{R^{k+1}}.
	\]
\end{theorem}

The existence of such a subsequence follows
from Theorem~\ref{L2}.
Along this subsequence, each individual factor
attains its extremal limit on $E$.
The fixed-row comparison and the Saff--Gonchar
theorem then yield the corresponding limits
on $E_R$. The partial-product and quotient
limits follow immediately.
The proof follows the argument
in~\cite[Section~4.7]{Pro2026},
where the details can be found.

The remainder of the paper is organized as follows.
In Section~4, we prove the comparison result,
Theorem~\ref{L6}.
In Section~5, we establish the ordinary limits
in Theorem~\ref{L12}, which are needed to prove
the main result of the paper,
Theorem~\ref{L2}, in Section~6.

\section{Proof of Theorem~\ref{L6}}

For each $R>1$, denote the Green level set by
\[
L_R:=\partial G_R
=\{z\in\Omega:g_\Omega(z,\infty)=\log R\}.
\]

To prove Theorem~\ref{L6}, it suffices
to establish the upper estimate \eqref{L7},
since the equality in \eqref{L8}
has already been deduced from this estimate
and the Saff--Gonchar theorem in Section~3.
We establish this estimate by polynomial
interpolation at the Fekete points of $E$.
Starting from a best rational
approximant $r=p/q$ on $E_R$, we interpolate
$qf-p$ at the Fekete points of $E$
to construct a rational function $r_1$
of the same class $\mathcal R_{n,m}$.
The approximation error of $r_1$ on $E$
provides an upper bound for
$\rho_{n,m}(f;E)$. 
We estimate this error using the interpolation
formula on a Green-level curve $L_s\subset E_R$,
obtaining a bound in terms of
$\rho_{n,m}(f;E_R)$, with the remaining factors
controlled by an elementary estimate for
the denominator $q$ and the standard
asymptotics of Fekete polynomials.

\begin{proof}
	
	Fix $1<s<R$ such that $\log s$ is a regular value of
	$g_\Omega(\cdot,\infty)$.
	The critical values of the Green function
	are locally finite, and $L_s$ is a finite union of analytic
	Jordan curves. We assume that $L_s$ is positively oriented
	with respect to $G_s$.
	
	Let $r=p/q\in\Rcal_{n,m}$
	be a best rational approximant to $f$ on $E_R$,
	where $p$ and $q$ have no common zeros. Then
	\[
	\|f-r\|_{E_R}=\rho_{n,m}(f;E_R).
	\]
	Since $f$ is holomorphic in a neighborhood of $E_R$
	and the approximation error is finite,
	$q$ has no zeros on $E_R$.
	
We begin with an elementary estimate for the
polynomial denominator $q$.

Set
\[
d=\dist(E,L_R)>0,
\qquad D=\diam E_R,
\]
and define
\[
C_{m,R}:=\left(1+\frac{D}{d}\right)^m.
\]

Since $q$ has degree at most $m$ and has no zeros
on $E_R$, we have
\begin{equation}\label{L11}
	\frac{\norm{q}{E_R}}{\min_E|q|}
	\leq C_{m,R}.
\end{equation}
Indeed, write
\[
q(z)=c\prod_{j=1}^{\ell}(z-b_j),
\qquad \ell\leq m,
\]
where $b_j\notin E_R$. For $x\in L_R$ and $y\in E$,
we have $|y-b_j|\geq d$, and hence
\[
\frac{|x-b_j|}{|y-b_j|}
\leq 1+\frac{|x-y|}{|y-b_j|}
\leq 1+\frac{D}{d}.
\]
Multiplying over all zeros and taking the maximum over $x\in L_R$ and the minimum over $y\in E$, we obtain
\[
\frac{\max_{L_R}|q|}{\min_E|q|}\le C_{m,R}.
\]
By the maximum modulus principle,
\[
\frac{\norm{q}{E_R}}{\min_E|q|}
=
\frac{\max_{L_R}|q|}{\min_E|q|}
\le C_{m,R},
\]
which proves \eqref{L11}.
For constant polynomials the
assertion is immediate. 

Let $\omega_{n+1}$ be the Fekete polynomial for $E$
of degree $n+1$.
Recall that its zeros belong to $E$.
Set $H=qf-p$, and let $P$ be the polynomial
of degree at most $n$ interpolating $H$
at the zeros of $\omega_{n+1}$. The interpolation formula, applied over the
contour $L_s$, gives
\[
H(z)-P(z)=\frac{1}{2\pi i}\int_{L_s}
\frac{\omega_{n+1}(z)}{\omega_{n+1}(t)}
\frac{H(t)}{t-z}\,dt,\qquad z\in G_s.
\]
Since $L_s\subset E_R$ and $H=q(f-r)$,
\[
\norm{H}{L_s}\le \norm{H}{E_R}
\le \norm{q}{E_R}\norm{f-r}{E_R}.
\]
Thus $r_1=(p+P)/q\in\Rcal_{n,m}$, and
\[
\norm{f-r_1}{E}
\le C_{E,s}\frac{\norm{\omega_{n+1}}{E}}
{\min_{L_s}|\omega_{n+1}|}
\frac{\norm{q}{E_R}}{\min_E|q|}
\norm{f-r}{E_R}, 
\]
where
\[
C_{E,s}:=
\frac{\operatorname{length}(L_s)}
{2\pi\dist(E,L_s)}
\]
depends only on $E$ and $s$ and is independent of $n$,
$m$, and the rational approximant $r$.
Using \eqref{L11} and the fact that $r$ is
a best rational approximant on $E_R$,
we obtain
\[
\rho_{n,m}(f;E)
\le C_{E,s}C_{m,R}
\frac{\norm{\omega_{n+1}}{E}}
{\min_{L_s}|\omega_{n+1}|}
\rho_{n,m}(f;E_R).
\]
The standard Fekete asymptotics (see
\cite{ST,Saff2010}) imply
\[
\lim_{n\to\infty}\norm{\omega_{n+1}}{E}^{1/n}
=\capE,
\qquad
\lim_{n\to\infty}
\left(\min_{L_s}|\omega_{n+1}|\right)^{1/n}
=s\capE,
\]
where $\capE$ denotes the logarithmic
capacity of $E$.
Dividing by $\rho_{n,m}(f;E_R)$, we obtain
\[
\frac{\rho_{n,m}(f;E)}
{\rho_{n,m}(f;E_R)}
\le C_{E,s}C_{m,R}
\frac{\norm{\omega_{n+1}}{E}}
{\min_{L_s}|\omega_{n+1}|}.
\]
Taking $n$th roots and upper limits gives
\[
\limsup_{n\to\infty}
\left(
\frac{\rho_{n,m}(f;E)}
{\rho_{n,m}(f;E_R)}
\right)^{1/n}
\le\frac1s.
\]
Letting $s\uparrow R$ through regular values
proves \eqref{L7}.
\end{proof}

\section{Existence of ordinary limits for rational approximation errors}

The proof of the Hadamard-type product formula
requires an extension of the Saff--Gonchar
asymptotic formula for rows whose meromorphy
radii are strictly smaller than $R_m$.
We establish that, whenever
$R_k<R_m$, the upper limit
\[
\limsup_{n\to\infty}
\rho_{n,k}(f;E)^{1/n}=\frac{1}{R_k}
\]
is in fact an ordinary limit. These limits
provide the essential ingredient for the
proof of the product formula in Section~6.

Saff~\cite{Saff1969b} obtained ordinary-limit
results for rows of the Walsh table under
additional assumptions on the geometry of
$E$ and the distribution of poles. 
Our proof is based on ideas developed in
Saff's papers~\cite{Saff1969a,Saff1969b},
as discussed in Section~2.
We use best rational approximants
from different rows of the Walsh table.
One of these rows is the first whose
meromorphy radius equals $R_m$,
so that Gonchar's analogue of the
de Montessus de Ballore theorem applies
to its approximants. 
An important step in the argument is the
application of the Bernstein--Walsh inequality.
Together with Rouch\'e's theorem and
pole counting, this leads to a contradiction
and yields the required limits.

We now recall Gonchar's analogue~\cite{G1974}
of the de Montessus de Ballore theorem~\cite{Mont1902}
for best rational approximants,
which extends the classical convergence theorem
for rows of the Pad\'e table
(see also~\cite{Gon1975a,Gon1975b,Gon1981}).

\newpage

\noindent{\bf Theorem. Gonchar's analogue of the
	de Montessus de Ballore theorem}
{\it Let $m\ge 1$ and $R>1$. Suppose that $f$
	admits a meromorphic continuation to $G_R$
	with exactly $m$ poles, counted with multiplicity.
	Let $r_{n,m}$ be best rational approximants
	of type $(n,m)$ to $f$ on $E$.
	Then the sequence $\{r_{n,m}\}$ converges
	uniformly to $f$ on compact subsets of $G_R$
	containing no poles of $f$.
	Moreover, for sufficiently large $n$,
	the approximants $r_{n,m}$ have exactly $m$
	poles, and each pole of $f$ in $G_R$
	of multiplicity $\ell$ attracts exactly
	$\ell$ poles of $r_{n,m}$ as $n\to\infty$. 
	}

We remark that the monic denominators converge
uniformly on compact subsets of $\mathbb C$
to the monic polynomial whose zeros are
the poles of $f$ in $G_R$, counted with
multiplicity.

We now establish the existence of the ordinary
limits in the following theorem.

\begin{theorem}\label{L12}
	Let $m\geq 1$ and suppose that $R_0<R_m$.
	Then, for every integer $k$ with
	$0\leq k\leq m-1$ and $R_k<R_m$,
	\[
	\lim_{n\to\infty}
	\rho_{n,k}(f;E)^{1/n}
	=\frac{1}{R_k}.
	\]
\end{theorem}

\begin{proof}
	Define
	\[
	p:=\max\{j\in\{0,1,\ldots,m-1\}:R_j<R_m\},
	\]
	so that
	\[
	R_0\leq\cdots\leq R_p
	<R_{p+1}=\cdots=R_m.
	\]
	Fix $k\in\{0,\ldots,p\}$. By the
	Saff--Gonchar theorem,
	\[
	\limsup_{n\to\infty}\rho_{n,k}(f;E)^{1/n}
	=\frac{1}{R_k}.
	\]
	Suppose, to the contrary, that
	\[
	\liminf_{n\to\infty}
	\rho_{n,k}(f;E)^{1/n}
	<\frac{1}{R_k}.
	\]
	Then there exist $R$ with $R_k<R<R_{p+1}=R_m$
	and an infinite subsequence $\Lambda$ such that
	\begin{equation}\label{L20}
		\rho_{n,k}(f;E)\le R^{-n},
		\qquad n\in\Lambda.
	\end{equation}
	Since $R_p<R_{p+1}$, the function $f$
	has exactly $p+1$ poles, counted with
	multiplicities, in every $G_s$ with
	$R_p<s<R_{p+1}$. Let $\omega$ be the monic
	polynomial whose zeros are these poles,
	counted with multiplicities.
		
	Write the best rational approximants
	in the $k$th and $(p+1)$st rows as
	\[
	r_{n,k}=\frac{p_n}{q_n},
	\qquad
	r_{n,p+1}=\frac{P_n}{Q_n},
	\]
	where the numerator and denominator in each
	fraction are relatively prime,
	$\|q_n\|_E=1$, and $Q_n$ is monic.

	Since $\deg q_n\le k$ and $\|q_n\|_E=1$,
	the polynomials $q_n$ form a bounded family
	in the finite-dimensional space of
	polynomials of degree at most $k$.
	Passing to a subsequence $\Lambda_1\subset\Lambda$,
	we may assume that
	\[
	q_n\longrightarrow q
	\quad\text{uniformly on compact subsets of }\mathbb C,
	\qquad n\to\infty,\quad n\in\Lambda_1,
	\]
	where
	\[
	\deg q\le k,\qquad \|q\|_E=1.
	\]
	In particular, $q\not\equiv0$.

	By Gonchar's theorem,
	\[
	Q_n\longrightarrow\omega,
	\qquad
	P_n\longrightarrow F:=f\omega
	\qquad n\to\infty,
	\]
	uniformly on compact subsets of $\mathbb C$
	and on compact subsets of $G_{R_{p+1}}$, respectively.
	Here $G_\infty:=\mathbb C$ if
	$R_{p+1}=\infty$.
	The function $F$ is holomorphic in
	$G_{R_{p+1}}$ and does not vanish
	at the poles of $f$.

	Since $R<R_m$, the Saff--Gonchar theorem gives, for some constant $C_0>0$
	independent of $n$,
	\begin{equation}\label{L21}
		\rho_{n,p+1}(f;E)\le C_0R^{-n},
		\qquad n\ge0.
	\end{equation}
	
	Set
	\[
	H_n=P_nq_n-p_nQ_n.
	\]
	Then $\deg H_n\le n+p+1$.
By the triangle inequality and $\|q_n\|_E=1$,
we have
\[
\|H_n\|_E
\le \|Q_n\|_E
\bigl(\rho_{n,k}(f;E)
+\rho_{n,p+1}(f;E)\bigr).
\]
Since $Q_n\to\omega$ uniformly on $E$ as
$n\to\infty$, the norms $\|Q_n\|_E$
are bounded.

Consequently, by \eqref{L20} and \eqref{L21},
there exists a constant $C>0$, independent
of $n$, such that
\[
\|H_n\|_E\le CR^{-n},
\qquad n\in\Lambda_1.
\]
Choose $R'$ with $R_k<R'<R<R_{p+1}$.
By the Bernstein--Walsh inequality,
\[
\|H_n\|_{E_{R'}}
\le C(R')^{p+1}
\left(\frac{R'}R\right)^n,
\qquad n\in\Lambda_1.
\]
Since $R'<R$, it follows that
\[
\|H_n\|_{E_{R'}}\longrightarrow0
\quad\text{as }n\to\infty,\quad n\in\Lambda_1.
\]
The inequalities $R_k<R'<R_{p+1}$ imply that
$f$ is meromorphic in $G_{R'}$.
Thus, $G_{R'}$ contains poles of $f$
of total multiplicity $N\ge k+1$.
Indeed, if $f$ had at most $k$ poles there, then 
$R_k\ge R'$, a contradiction.

Since $R'<R_{p+1}$, the convergence
$P_n\to F$ also holds uniformly on compact
subsets of $G_{R'}$.
Using $p_nQ_n=P_nq_n-H_n$, we obtain
\[
p_nQ_n\longrightarrow Fq
\]
uniformly on compact subsets of $G_{R'}$ as
$n\to\infty$, $n\in\Lambda_1$.

By Gonchar's theorem, each pole of $f$
in $G_{R'}$ of multiplicity $\ell$
attracts exactly $\ell$ zeros of $Q_n$.
Choose sufficiently small disjoint disks
around these poles, with closures contained
in $G_{R'}$, such that $F$ has no zeros
in their closures.

Since $p_nQ_n\to Fq$ uniformly on compact
subsets of $G_{R'}$, Rouch\'e's theorem
implies that $Fq$ has at least $\ell$
zeros in each disk corresponding to a pole
of multiplicity $\ell$.
Since $F$ does not vanish at the poles
and the disks may be chosen arbitrarily small,
$q$ must have a zero at each pole
of multiplicity at least equal to that
of the pole. Consequently,
\[
k\ge\deg q\ge N\ge k+1,
\]
a contradiction.
Therefore,
\[
\liminf_{n\to\infty}
\rho_{n,k}(f;E)^{1/n}\ge\frac1{R_k}.
\]
Together with the Saff--Gonchar theorem,
this proves the theorem.
\end{proof}

    \section{Proof of Theorem~\ref{L2}}
    
    The key ingredient in the proof is
    Theorem~\ref{L12}, which establishes
    ordinary limits for the rows whose
    meromorphy radii are strictly smaller
    than $R_m$. For the remaining rows,
    the inclusion
    \[
    \mathcal R_{n-m+k,k}\subseteq\mathcal R_{n,m},
    \qquad 0\leq k\leq m,
    \]
    provides lower bounds in terms of the
    approximation error in the $m$th row.
    Combining these observations with the
    Saff--Gonchar theorem yields the required
    lower bound for the product.

    \begin{proof}
    	The upper estimate in \eqref{L3}
    	has already been established in
    	\eqref{L1}.
    	It remains to prove the reverse inequality.
    	We distinguish three cases.
    	
    	\medskip
    	\noindent\textbf{Case 1: $R_m=\infty$.}

    In this case, the right-hand side of
    	\eqref{L1} is zero.
    	Since the approximation errors are
    	nonnegative, it follows that
    	\[
    	\lim_{n\to\infty}
    	\left(
    	\prod_{k=0}^{m}
    	\rho_{n-m+k,k}(f;E)
    	\right)^{1/n}=0.
    	\]
    	Thus, \eqref{L3} holds.

    	\medskip
    	\noindent\textbf{Case 2: $R_m<\infty$ and
    		$R_0=R_m$.}
    	
    	By the monotonicity of the meromorphy radii,
    	\[
    	1<R_0\le R_1\le\cdots\le R_m,
    	\]
    	the equality $R_0=R_m$ implies that
    	\[
    	R_0=R_1=\cdots=R_m.
    	\]
    	In particular,
    	\[
    	R_0R_1\cdots R_m=R_m^{m+1}.
    	\]
    	Since
    	\[
    	\mathcal R_{n-m+k,k}
    	\subseteq\mathcal R_{n,m},
    	\qquad 0\le k\le m,
    	\]
    	we have
    	\begin{equation}\label{L14}
    		\rho_{n-m+k,k}(f;E)
    		\ge \rho_{n,m}(f;E).
    	\end{equation}
    	Therefore,
    	\[
    	\prod_{k=0}^{m}
    	\rho_{n-m+k,k}(f;E)
    	\ge
    	\rho_{n,m}(f;E)^{m+1}.
    	\]
    	Taking $n$th roots and upper limits,
    	and using the Saff--Gonchar theorem,
    	we obtain
    	\[
    	\begin{aligned}
    		&\limsup_{n\to\infty}
    		\left(
    		\prod_{k=0}^{m}
    		\rho_{n-m+k,k}(f;E)
    		\right)^{1/n}
    		\\
    		&\qquad\ge
    		\limsup_{n\to\infty}
    		\rho_{n,m}(f;E)^{(m+1)/n}
    		\\
    		&\qquad=
    		\frac{1}{R_m^{m+1}}
    		=
    		\frac{1}{R_0R_1\cdots R_m}.
    	\end{aligned}
    	\]
    	Together with \eqref{L1},
    	this proves \eqref{L3} in the present case.
    	
    	\medskip
    	\noindent\textbf{Case 3: $R_0<R_m<\infty$.}
    	
    	Define
    	\[
    	p:=\max\{j\in\{0,1,\ldots,m-1\}:R_j<R_m\}.
    	\]
    	By the monotonicity of the meromorphy radii
    	and the maximality of $p$,
    	\begin{equation}\label{L15}
    		1<R_0\le R_1\le\cdots\le R_p
    		<R_{p+1}=\cdots=R_m.
    	\end{equation}
    	By Theorem~\ref{L12}, for every
    	$k=0,1,\ldots,p$, the full limit exists and
    	satisfies
    	\[
    	\lim_{n\to\infty}
    	\rho_{n,k}(f;E)^{1/n}
    	=\frac{1}{R_k}.
    	\]
    	Since $m$ and $k$ are fixed, replacing $n$
    	by $n-m+k$ does not change the exponential
    	rate. Hence
    	\[                                         	\lim_{n\to\infty}
    		\rho_{n-m+k,k}(f;E)^{1/n}
    		=\frac{1}{R_k},
    		\qquad 0\le k\le p.
    	\]
    	Multiplying these limits, we obtain
    	\begin{equation}\label{L16}
    		\lim_{n\to\infty}
    		\left(
    		\prod_{k=0}^{p}
    		\rho_{n-m+k,k}(f;E)
    		\right)^{1/n}
    		=
    		\frac{1}{R_0R_1\cdots R_p}.
    	\end{equation}
    	On the other hand, by
    	\eqref{L14},
    	\[
    	\rho_{n-m+k,k}(f;E)
    	\ge \rho_{n,m}(f;E),
    	\qquad p+1\le k\le m.
    	\]
    	We now estimate the upper limit of the $n$th roots of the full product from below.
    	Combining these inequalities with
    	\eqref{L16} and the Saff--Gonchar theorem,
    	we obtain
   
    \[
    \begin{aligned}
    	&\limsup_{n\to\infty}
    	\left(
    	\prod_{k=0}^{m}
    	\rho_{n-m+k,k}(f;E)
    	\right)^{1/n}
    	\\
    	&\quad\ge
    	\frac{1}{R_0R_1\cdots R_p}
    	\limsup_{n\to\infty}
    	\rho_{n,m}(f;E)^{(m-p)/n}
    	\\
    	&\quad=
    	\frac{1}{R_0R_1\cdots R_pR_m^{m-p}}
    	\\
    	&\quad=
    	\frac{1}{R_0R_1\cdots R_m},
    \end{aligned}
    \]
    where the last equality follows from
    \eqref{L15}.
    
    Combining this lower bound with
    \eqref{L1}, we obtain
    \eqref{L3} in the remaining case.
    The proof is complete.
    \end{proof}

 \bigskip
 \noindent Department of Mathematics and Statistics,\\
 University of South Alabama,\\
 Mobile, AL 36688-0002, USA.
 
 \noindent E-mail: \texttt{prokhoro@southalabama.edu}

\end{document}